\documentclass[12
pt,twoside,a4paper]{amsart}
\usepackage{amssymb}
\usepackage[latin1]{inputenc}
\usepackage[T1]{fontenc}
\usepackage{times}

\def\hpq0{h^{p,q}_{\leq 0}}
\def\Hpq0{\H_{\leq 0}^{p,q}}
\def\gr{\otimes}

\def\dbar{\bar\partial}
\def\ddbar{\partial\dbar}

\def\C{{\mathbb C}}

\def\H{{\mathcal H}}

\def\be{\begin{equation}}
\def\ee{\end{equation}}

\newtheorem{thm}{Theorem}[section]

\theoremstyle{definition}

\theoremstyle{remark}

\newtheorem{preremark}{Remark}
\newtheorem{preex}{Example}

\numberwithin{equation}{section}

\begin{document}

\title[]
{Curvature of vector bundles and subharmonicity of Bergman kernels.}

\author[]{ Bo Berndtsson}

\address{B Berndtsson :Department of Mathematics\\Chalmers University
  of Technology 
  and the University of G\"oteborg\\S-412 96 G\"OTEBORG\\SWEDEN,\\}

\email{ bob@math.chalmers.se}

\begin{abstract}
{In a previous paper, \cite{Berndtsson}, we have studied a property of
  subharmonic 
  dependence on a parameter of Bergman kernels for a family of weighted
  $L^2$-spaces of holomorphic 
  functions. Here we prove a result on the
  curvature of a vector bundle defined by this family of $L^2$-spaces itself,
  which has the earlier results on Bergman kernels as a corollary.
  Applying the same arguments to spaces of holomorphic sections to
  line bundles over a locally trivial fibration we also prove that if
  a holomorphic vector bundle, $V$, over a complex manifold is ample
  in the sense of 
  Hartshorne, then $V\gr\det V$ has an Hermitian metric with curvature
  strictly positive in the sense of Nakano.}
\end{abstract}

\bigskip

\maketitle

\section{Introduction}

Let us consider a domain $D=U\times\Omega$  and a 
function $\phi$ plurisubharmonic in $D$. We also assume
for simplicity that $\phi$ is smooth up to the boundary and  strictly
plurisubharmonic in $D$. Then, for each $t$ in $U$,
$\phi^t(\cdot):=\phi(t,\cdot)$ is plurisubharmonic in $\Omega$ and we denote
by $A^2_t$ the Bergman spaces of holomorphic functions in $\Omega$
with norm
$$
\|h\|^2=\|h\|^2_t=\int_\Omega |h|^2e^{-\phi^t}. 
$$
The spaces $A^2_t$ are then all equal as vector spaces but have norms
that vary with $t$. The - infinite rank - vector bundle $E$ over $U$
with fiber $E_t=A^2_t$ is therefore trivial as a bundle but is
equipped with a nontrivial metric. The main result of this paper is
the following theorem. 
\begin{thm} The hermitian bundle $(E, \|\cdot\|_t)$ is strictly positive
  in the sense of Nakano.
\end{thm}

Of the two main differential geometric notions of positivity (see
section 2, where these matters will be reviewed in the slightly non
standard setting of bundles of infinite rank), positivity in the sense
of Nakano is the stronger one and 
 implies the  weaker property of positivity in the sense of
Griffiths. On the other hand the Griffiths notion of positivity has
nicer functorial properties and implies in particular that the dual
bundle is negative (in the sense of Griffiths). This latter property
is in turn equivalent to the condition that if $\xi$ is  any
nonvanishing local
holomorphic section to the dual bundle, then  the function
$$
\log\|\xi\|^2
$$
is strictly plurisubharmonic. In our case, we can obtain such
holomorphic sections to the dual bundle from point evaluations. More
precisely, let $f$ be a holomorphic map from $U$ to $\Omega$ and
define $\xi_t$ by is action on a local section to $E$ 
$$
\langle \xi_t, h_t\rangle= h_t(f(t)).
$$
Since the right hand side here is a holomorphic function of $t$, $\xi$
is indeed a holomorphic section to $E^*$. The norm of $\xi$ at a point
is given by
$$
\|\xi_t\|^2=\sup_{\|h_t\|\leq 1}| h_t(f(t))|=K_t(f(t),f(t)),
$$
where $K_t(z,z)$ is the Bergman kernel function for $A^2_t$. It
therefore follows from Theorem 1.1 that $K_t(z,z)$ is plurisubharmonic
in $D$, which is essentially the main result of \cite{Berndtsson}. 

The proof of Theorem 1.1 is based on a formula of Griffiths, see e g
\cite{Griff-Harr},  for the
curvature of a  subbundle, $E$, of a holomorphic bundle, $F$. If
$\Theta^E$ and  $\Theta^F$ denote the respective curvature operators
we have, for any sections    $u$ and $v$ of the subbundle $E$ that 
$$
(\Theta^F_{j k}u,v)=(\Theta^E_{j k}u,v)
+(\pi_{E^\bot}D^F_{j}u,\pi_{E^\bot}D^F_{k}v) .
$$
In this formula, the last term equals the {\it second fundamental
  form} of $E$ acting on $u$ and $v$ , and is here expressed by 
 the Chern connection of $F$ acting on sections of $E$,  projected on the
orthogonal 
complement of $E$ in $F$. 

We shall apply this formula with $E$ being the bundle introduced above
and $F$ the bundle of all $L^2$ functions on $\Omega$ equipped with
the same norm. The curvature of $F$ is  easily seen to be the
operator of multiplication with the Hessian of $\phi$ with respect to
$t$. This is therefore a positive operator, and to prove Theorem 1.1, we
must control the second fundamental form with this operator. For this
we note that the second fundamental 
form is, by definition, given by  elements in the orthogonal complement of
$A^2_t$. These elements are  the $L^2$-minimal solutions of  certain
$\dbar$-equations. An estimate for these $L^2$-minimal solutions is
furnished by Hörmander's $L^2$-estimate for the $\dbar$-equation. The
basic idea behind the proof comes from a proof of {\it Prekopa's theorem}
by Brascamp and Lieb, and is explained more closely in \cite{Berndtsson}.

There is also a natural analog of Theorem 1.1 for locally trivial
fibrations. We  consider a complex manifold $X$ which is fibered
over 
another complex manifold $Y$. We then have a holomorphic map, $p$, from
$X$ to $Y$ with surjetive differential, and all the fibers
$X_t=p^{-1}(t)$ are diffeomorphic. We shall even assume that the fibration
is locally trivial holomorphically, so that every point in the base
has a neighbourhood, $U$, such that $p^{-1}(U)$ is biholomorphic to
$U\times Z$, where $Z$ is a fixed  complex manifold. Moreover, under
this biholomorphism, the projection $p$ goes over into the natural
projection from $U\times Z$ to $U$.
Let $L$ be a positive line bundle over $X$ and assume that under the
local trivializations discussed above, $L$ restricted to $p^{-1}(U)$
is isomorphic to  the pullback of a line bundle on $Z$ under the
natural projection from $U\times Z$ to $Z$. For each $t$ in $Y$ we can
now consider the space
$$
\tilde E_t=\Gamma(X_t,L|X_t)
$$
of global holomorphic sections to $L$ over $X_t$. By the local
triviality of the fibration, and the assumption on $L$, the vector spaces
$\tilde E_t$ define a holomorphic vector bundle over $Y$. We would now like
to define a hermitian norm on the bundle $\tilde E$ by taking the $L^2$-norm
over each fiber, but we can not do so directly since we have no
canonically defined measure on the fibers to integrate against. We therefore
consider instead the bundle $E$ with fibers
$$
E_t=\Gamma(X_t,L|X_t\gr K_{X_t}),
$$
where $K_{X_t}$ is the canonical bundle of each fiber. Elements of
$E_t$ can be naturally integrated over the fiber and we obtain in this
way a metric, $\|\cdot\|$ on $E$ in complete analogy with the plane
case. We then get the same conclusion as before:
\begin{thm} $(E,\|\cdot\|)$ is positive in the sense of Nakano.
\end{thm}
One example of this situation that arises naturally is obtained if we
start with a (finite rank) holomorphic vector bundle $V$ over $Y$ and
let $\mathbb P(V)$ be the associated bundle of projective spaces {\it of
the dual bundle} $V^*$. This is then clearly a locally trivial
holomorphic fibration as before and a line bundle $L$ satisfying the
conditions we have discussed is obtained by taking
$$
L=O_{\mathbb P(V)}(1),
$$
the hyperplane section bundle over each fiber. The global holomorphic
sections of 
this bundle over each fiber are now the linear forms on $V^*$, i e the
elements of $V$. In other words, $\tilde E$ is isomorphic to $V$. As
explained above, we are not able to produce a metric on $\tilde E$ by
integrating over the fibers, so instead we take as $L$
$$
L=O_{\mathbb P(V)}(r+1)
$$
(with $r$ being the rank of $V$) and define $E$ as before
$$
E=\Gamma(X_t,L|X_t\gr K_{X_t}).
$$
One can then verify that $E$ is isomorphic to $V\gr\det V$. The
condition that $L$ is positive is now equivalent to $O_{\mathbb P(V)}(1)$ being
positive which is the same as saying that $V$ is {\it ample} in the
sense of Hartshorne, \cite{Harts}. We therefore obtain the following result
 as a corollary of Theorem 1.2.

\begin{thm} Let $V$ be a (finite rank) holomorphic vector bundle over
  a complex manifold which is ample in the sense of Hartshorne. Then
  $V\gr\det V$ has a smooth hermitian metric which is strictly
  positive in the sense of Nakano. 
\end{thm}
It is a well known conjecture of Griffiths, \cite{Griffiths}, that an ample
vector bundle is positive in the sense of Griffiths.
Theorem 1.3 would follow from this conjecture, since by a a theorem of
Demailly, \cite{Demailly 2}, $V\gr\det V$ is Nakano positive if $V$ itself is
Griffiths positive. It seems however that not so much is known about
Griffiths' conjecture in general, except that it does hold when $Y$ is a
compact curve (see \cite{Umeu}, \cite{Camp}).

After this mansucript was completed I received a preprint
\cite{Mourougane-Takayama}. They prove that $V\gr\det V$ is positive in
the sense of Griffiths, assuming the base manifold is projective. The
method of proof seems quite different.  
Finally, I would like to thank Sebastien Boucksom for pointing out the
relation between Theorem 1.1 and the Griffiths conjecture.

\section{Curvature of finite and infinite rank bundles}

Let $E$ be a holomorphic vector bundle with a hermitian metric over a
complex manifold $Y$. By 
definition this means that there is a holomorphic projection map $p$
from $E$ to $Y$ and that every point in $Y$ has a neighbourhood $U$
such that $p^{-1}(U)$ is isomorphic to $U\times W$, where $W$ is a
vector space equipped with a smoothly varying hermitian metric. In our
applications it is important to be able to allow this vector space to
have infinite dimension, in which case we assume that the metrics are
also complete, so that the fibers are Hilbert spaces. 

Let $t=(t_1,...t_m)$ be a system of local coordinates on $Y$. The
Chern connection, $D_{t_j}$ is now given by a collection of
differential operators acting on smooth sections to $U\times W$ and
satisfying
$$
\partial_{t_j}( u,v)=(D_{t_j}u,v) +( u,\dbar_{t_j}v),
$$ 
with $\partial_{t_j}=\partial/\partial_{t_j}$ and
$\dbar_{t_j}=\partial/\partial\bar t_j$. The curvature of the Chern 
connection is a $(1,1)$-form of operators 
$$
\Theta=\sum\Theta_{j k}dt_j\wedge d\bar t_k,
$$
where the coefficients $\Theta_{j k}$ are densily defined operators on
$W$. By definition these coefficients are the commutators
$$
\Theta_{j k}= [ D_{t_j},\dbar_{t_k}].
$$
The vector bundle is said to be positive in the sense of Griffiths if
for any section $u$  to $W$ and any vector $v$ in $\C^m$
$$
\sum(\Theta_{j k}u,u)v_j\bar v_k\geq \delta\|u\|^2|v|^2
$$
for som positive $\delta$. $E$ is said to be positive in the sense of
Nakano if for any $m$-tuple $(u_1,...u_m)$ of sections to $W$
$$
\sum(\Theta_{j k}u_j,u_k)\geq \delta\sum\|u_j\|^2
$$
Taking $u_j=u v_j$ we see that Nakano positivity implies positivity in
the sense of Griffiths. 

The dual bundle of $E$ is the vector bundle $E^*$ whose fiber at a
point $t$ in $Y$ is the Hilbert space dual of $E_t$. There is
therefore a natural  antilinear isometry
between $E$ and $E^*$, which we will denote by  $J$. If $u$ is a local
section to $E$, $\xi$ is 
a local section to $E^*$, and $\langle\cdot,\cdot\rangle$ denotes the
pairing between $E^*$ and $E$ we have
$$
\langle\xi,u\rangle=(u,J\xi).
$$
Under the natural holomorphic structure on $E^*$ we then have
$$
\dbar_{t_j}\xi=J^{-1}D_{t_j}J\xi,
$$
and the Chern connection on $E^*$ is given by
$$
D^*_{t_j}\xi=J^{-1}\dbar_{t_j}J\xi.
$$
It follows that
$$
\dbar_{t_j}\langle\xi,u\rangle=
\langle\dbar_{t_j}\xi,u\rangle +\langle\xi,
\dbar_{t_j}u\rangle,
$$
and 
$$
\partial_{t_j}\langle\xi,u\rangle=
\langle D^*_{t_j}\xi,u\rangle +\langle\xi,D_{t_j}u\rangle,
$$
and hence 
$$
0=\left[\partial_{t_j},\dbar_{t_j}\right]\langle\xi,u\rangle =
\langle\Theta^*_{j k}\xi,u\rangle +\langle\xi,\Theta_{j k}u\rangle,
$$ 
if we let $\Theta^*$ be the curvature of $E^*$.
If $\xi_j$ is an $r$-tuple of sections to $E^*$, and $u_j=J\xi_j$,  we
thus see that 
$$
\sum(\Theta^*_{j k}\xi_j,\xi_k)=-\sum(\Theta_{j k}u_k,u_j).
$$
Notice that the order between $u_k$ and $u_j$ in the right hand side is
opposite to the order between the $\xi$s in the left hand
side. Therefore $E^*$ is negative in the sense of Griffiths iff $E$ is
positive in the sense of Griffiths, but we can not draw the same
conclusion in the the case of Nakano positivity.

If $u$ is a holomorphic section to $E$ we also find that
$$
\frac{\partial^2}{\partial t_j\partial\bar t_k}(u,u)=
(D_{t_j}u,D_{t_k}u) -(\Theta_{j k}u,u)
$$
and it follows after a short computation that $E$ is (strictly)
negative in the 
sense of Griffiths if and only if $\log\|u\|^2$ is (strictly)
plurisubharmonic for any nonvanishing holomorphic section $u$.

We next briefly recapitulate the Griffiths formula for the curvature
of a subbundle. Assume $E$ is a holomorphic subbundle of the bundle
$F$, and let $\pi$ be the fiberwise orthogonal projection from $F$ to
$E$. We also let $\pi_\bot$ be the orthogonal projection on the
orthogonal complement of $E$.
By the definition of Chern connection we have 
$$
D^E=\pi D^F.
$$
Let $\dbar_{t_j}\pi$ be defined by
\begin{equation}
\dbar_{t_j}(\pi u)=(\dbar_{t_j}\pi)u +\pi(\dbar_{t_j} u).
\end{equation}
Computing the commutators occuring in the  definition of curvature we see that
\begin{equation}
\Theta^E_{j k} u=-(\dbar_{t_k}\pi)D^F_{t_j} u+\pi\Theta^F_{j k} u,
\end{equation}
if $u$ is a section to $E$. 
By (2.1) $(\dbar\pi)v=0$ if $v$ is a  section to $E$,
so
\begin{equation}
(\dbar\pi)D^F u=(\dbar\pi)\pi_\bot D^F u.
\end{equation}
 Since
$\pi \pi_\bot=0$ it also follows that
$$
(\dbar\pi)\pi_\bot D^F u=-\pi\dbar(\pi_\bot D^Eu),
$$
so  if $v$ is also a section to $E$,
$$
((\dbar_{t_k}\pi)D^F_{t_j} u,v)=-(\dbar_{t_k}(\pi_\bot D^F_{t_j}u),v)=
$$
$$
=
((\pi_\bot
D^F_{t_j}u),D^F_{t_k}v)=(\pi_\bot(D^F_{t_j}u),\pi_\bot(D^F_{t_k}v)) .
$$
Combining with (2.2) we finally get that if $u$ and $v$ are both
sections to $E$ 
then 
\begin{equation}
(\Theta^F_{j k}u,v)=(\pi_\bot(D^F_{t_j}u),\pi_\bot(D^F_{t_k}v))+
(\Theta^E_{j k}u,v),
\end{equation}
which is the starting point for the proof in the next section.

\section{ The proof of Theorem 1.1}

We consider the set up described before the statement of Theorem 1.1
in the introduction. Thus $E$ is the vector bundle over $U$ whose
fibers are the Bergman spaces $A^2_t$ equipped with the weighted $L^2$
metrics induced by  $L^2(\Omega,e^{-\phi^t})$. We also let $F$ be the
vector bundle with fiber  $L^2(\Omega,e^{-\phi^t})$, so that $E$ is a
trivial subbundle of the trivial bundle $F$ with a  metric
induced from a nontrivial metric on $F$. From the definition of the
Chern connection we see that
$$
D^F_{t_j}=\partial_{t_j}-\phi_j,
$$
where the last term in the right hand side should be interpreted as the operator of
multiplication by the (smooth) function
$-\phi_j=-\partial_{t_j}\phi^t$. (In the sequel we use the letters
$j, k$ for indices of the $t$-variables,and the letters
$\lambda,\mu$ for indices of the $z$-variables.) For the
curvature of $F$ we therefore get
$$
\Theta^F_{j k}=\phi_{j k},
$$
the operator of multiplication with the Hessian of $\phi$ with respect
to the $t$-variables. We shall now apply formula (2.4), so let $u_j$ be
 smooth sections to $E$. This means that $u_j$ are functions that
depend smoothly on $t$ and holomorphically on $z$. To verify the
positivity of $E$ in the sense of Nakano we need to estimate from
below the
curvature of $E$ acting on the $k$-tuple $u$,
$$
\sum(\Theta^E_{j k}u_j,u_k).
$$
By (2.4) this means that we need to estimate from above
$$
\sum (\pi_\bot(\phi_ju_j),\pi_\bot(\phi_k u_k))=\|\pi_\bot(\sum \phi_ju_j)\|^2.
$$
Put $w=\pi_\bot(\sum \phi_ju_j)$.  For fixed $t$, $w$ solves the
$\dbar_z$-equation 
$$
\dbar w=\sum u_j\phi_{j \lambda}d\bar z_\lambda,
$$
since the $u_j$s  are holomorphic in $z$. Moreover, since $w$ lies in
the orthogonal complement of $A^2$, $w$ is the minimal
solution to this equation. 

We shall next apply Hörmander's weighted $L^2$-estimates for the
$\dbar$-equation. The precise form of these estimates that we need
says that if $f$ is a $\dbar$-closed form in a psedudoconvex domain 
$\Omega$, and if $\psi$ is a smooth strictly plurisubharmonic weight
function, then the minimal solution $w$ to the equation $\dbar v=f$
satisfies
$$
\int_\Omega |w|^2 e^{-\psi}\leq \int_\Omega\sum\psi^{\lambda
  \mu}f_\lambda\bar f_\mu e^{-\psi},
$$
where $(\psi^{\lambda \mu})$ is the inverse of the complex Hessian of
$\psi$ (see \cite{Demailly 1}).

In our case this means that 
$$
\int_\Omega|w|^2 e^{-\phi^t}\leq\int_\Omega\sum \phi^{\lambda
  \mu}\phi_{j \lambda}u_j\overline{\phi_{k\mu}u_k} e^{-\phi^t}.
$$
Inserting this estimate in formula (2.4) together with the formula for
the curvature of $F$ we find
\begin{equation}
\sum(\Theta^E_{j k}u_j,u_k)\geq\int_\Omega \sum_{j k}\left(\phi_{j
    k}-\sum_{\lambda \mu}\phi^{\lambda \mu}\phi_{j
    \lambda}\bar\phi_{k\mu}\right) u_j\bar u_k e^{-\phi^t}.
\end{equation}
We claim that the expression
$$
D_{j k}=:\sum_{j k}\left(\phi_{j
    k}-\sum_{\lambda \mu}\phi^{\lambda \mu}\phi_{j
    \lambda}\bar\phi_{k\mu}\right),
$$
in the integrand is a positive definite matrix at any fixed point . By
a linear change
of variables in $t$ we may of course assume that the vector $u$ that
$D$ acts on equals $(1,0...0)$. Let $\Phi=i\ddbar \phi$ where the
$\ddbar$-operator acts on $t_1$ and the $z$-variables, the remaining
$t$-variables being fixed. Then
$$
\Phi=\Phi_{ 1 1} +i\alpha\wedge d\bar t_1+idt_1\wedge\bar\alpha +\Phi',
$$
where $\Phi_{1 1}$ is of bidegree $(1,1)$ in $t_1$, $\alpha$ is of
bidegree $(1,0)$  in $z$, and $\Phi'$ is of bidegree
$(1,1)$ in $z$. Then
$$
\Phi_{m+1}=\Phi^{m+1}/(m+1)!=
\Phi_{1,1}\wedge\Phi'_m-i\alpha\wedge\bar\alpha\wedge\Phi'_{m-1}\wedge
idt_1\wedge d\bar t_1.
$$
Both sides of this equation are forms of maximal degree that can be
written as certain coefficients multiplied by  the Euclidean volume form of
$\C^{m+1}$. The coefficient of the left hand side is the hessian of
$\phi$ with respect to $t_1$ and $z$ together. Similarily, the
coefficient of the first term on the right hand side is $\phi_{1 1}$
times the hessian of
$\phi$ with respect to the $z$-variables only. Finally, the
coefficient of the last term on the right hand side is the norm of the
$(0,1)$ form in $z$
$$
\dbar_z \partial_{t_1}\phi
$$
measured in the metric defined by $\Phi'$, multiplied by the volume
form of the same metric.
Dividing by the coefficient of $\Phi'_m$ we thus see that the matrix $D$
acting on a vector $u$ as above equals the hessian of $\phi$ with
respect to $t_1$ and $z$ divided by the hessian   of $\phi$ with respect to
the $z$-variables only. This expression is therefore positive so the
proof of Theorem 1.1 is complete. 

\section{Locally trivial fibrations} 

The proof of Theorem 1.2 is essentially the same as the proof in the
preceeding paragraph. As the statement is local we may assume that the
total space $X$ is a product $X=U\times Z$ where $U$ is open in $\C^k$
and $Z$ is a compact complex manifold of dimension $m$. Over $X$ we have a
holomorphic line bundle $L$ which is the pull back of a bundle on
$Z$. The line bundle $L$ is given a metric with strictly positive
curvature on $X$. Slightly abusively, we will denote by $K_Z$ the pull
back of the canonical bundle (i e the bundle of $(m,0)$-forms) on $Z$
to $X$ and also the restriction of this bundle to each fiber
$X_t=\{t\}\times Z$. On each fiber there is a natural pairing of sections
to $L\gr K_Z$ with values in the space of forms of maximal degree on
the fiber. Locally, for a decomposable section, $u=s\gr dz$ it is
given by
$$
[u,u]=c_k|s|^2 dz\wedge d\bar z,
$$
where $c_k$ is a constant of modulus 1 chosen so that the expression
is always nonnegative and  $|s|$ is the norm given by the metric on
$L$. For each $t$ in 
$U$ this induces a metric on the space of smooth section to $L\gr K_Z$
over $X_t$,
$$
\|u\|^2=\int_{X_t}[u,u].
$$
We let $F$ be the (infinite rank) vector bundle over $U$ whose
fiber over a point $t$ in $U$ is the $L^2$-space defined by this
metric. Similarily $E$ is the (finite rank)  vector bundle of holomorphic
sections. Note that, while $E$ naturally extends as a bundle over all
of $Y$, there seems to be no canonical way of extending the definition
of $F$ globally. Since our computation are local, this is however not
needed. 
Let $\theta$ be the connection form of the Chern connection on $L$
over $X$. Locally, in terms of a local trivialisation where the metric
is given by a weight function $\phi$, $\theta=-\partial\phi$. 
The Chern connection on the bundle $F$ is now
$$
D^F_{t_j}=\partial_{t_j}+\theta_{t_j},
$$
where $\theta_{t_j}$ is the coefficient of $dt_j$ in $\theta$. For the
curvature of $F$ we get
$$
\Theta^F_{j k}=c_{j k},
$$
the operator of multiplication with the $t$-part of the curvature of
$L$.
The proof of Theorem 1.2 now follows the same lines as the proof of
Theorem 1.1, using the Kodaira-Nakano-Hörmander estimate for line
bundles over compact manifolds, see \cite{Demailly 1}.  
\section{ Bundles of projective spaces}
Let $V$ be a  holomorphic line bundle of finite rank $r$ over a complex
manifold $Y$, and let $V^*$ be its dual bundle. We let $\mathbb P(V)$
be the fiber bundle over $Y$ whose fiber at each point $t$ of the base
is the projective space of lines in $V^*_t$, $\mathbb P(V^*_t)$. Then
$\mathbb P(V)$ is a holomorphically locally trivial fibration. There
is a naturally defined line bundle $O_{\mathbb P(V)}(1)$ over $\mathbb
P(V)$ whose restriction to any  fiber $\mathbb P(V^*_t)$ is the hyperplane
section bundle. One way to define this bundle is to first consider the
tautological line bundle
$O_{\mathbb P(V)}(-1)$. The total s   pace of this line bundle is just
the total space of $V^*$ with the zero section removed, and the
projection to  $\mathbb P(V)$ is the map that sends a nonzero point in
$V^*_t$ to its image in  $\mathbb P(V^*_t)$. The bundle 
$O_{\mathbb P(V)}(1)$ is then defined as the dual of $O_{\mathbb P(V)}(-1)$. 
The global holomorphic sections
of this bundle over any fiber are in one to one correspondence with
the linear forms on $V_t^*$, i e the elements of $V$. More generally, $
O_{\mathbb P(V)}(1)^l=O_{\mathbb P(V)}(l)$ has as global holomorphic
sections over each fiber the homogenuous polynomials on $V^*_t$ of
degree $l$, i e the elemets of the $l$:th symmetric power of $V$.
We shall apply Theorem 1.2 to the line bundles 
$$
L(l)=: O_{\mathbb P(V)}(l).
$$

Let  $E(l)$ be the vector bundle whose fiber over a point $t$ in $Y$ is
the space of global holomorphic sections of 
$L(l)\gr K_{\mathbb P(V^*_t)}$. If $l<r$ there is only the zero
section, so we assume 
from now on that $l$ is greater than or equal to $r$. 

We claim that
$$
E(r)=\det V,
$$
the determinant bundle of $V$. To see this, note that $L(r)\gr
K_{\mathbb P(V^*_t)}$ is trivial on each fiber, since the canonical
bundle of $(r-1)$-dimensional projective space is $O(-r)$. The space
of global sections is therefore one dimensional. A convenient  basis
element is
$$
\sum_1^r z_j \widehat{dz_j},
$$
if $z_j$ are coordinates on $V^*_t$. Here $\widehat{dz_j}$ is the
wedge product of all differentials $dz_k$ except $dz_j$ with a sign
chosen so that $dz_j\wedge\widehat{dz_j}= dz_1\wedge ...dz_r$. If we
make a linear change of 
coordinates on $V^*_t$, this basis element gets multiplied with the
determinant of the matrix giving the change of coordinates, so the
bundle of sections must transform as the determinant of $V$. Since
$$
L(r+1)\gr K_{\mathbb P(V^*_t)}=O_{\mathbb P(V)}(1)\gr L(r)\gr
K_{\mathbb P(V^*_t)},
$$
it also follows that
$$
E(r+1)=V\gr\det V.
$$ 
In the same way 
$$
E(r+m)= S^m(V)\gr\det V,
$$
where  $S^m(V)$ is the $m$th symmetric power of $V$. 

Let us now assume that $V$ is ample in the sense of Hartshorne, see
\cite{Harts}. By a theorem of Hartshorne, \cite{Harts}, $V$ is ample
if and only if 
 $L(1)$ is  ample, i e has a metric with
strictly positive curvature. Theorem 1.2 then implies that the
$L^2$-metric on each of the bundles $E(r+m)$
for $m\geq 0$ has curvature which is strictly positive in the sense of
Nakano, so we obtain:
\begin{thm} Let $V$ be a vector bundle (of finite rank) over a complex
  manifold. 
Assume $V$ is ample in the sense of Hartshorne. Then for any $m\geq 0$
the bundle
$$
S^m(V)\gr\det V
$$
has an hermitian metric with curvature which is (strictly) positive in
the sense of Nakano.
\end{thm}

\bigskip

\def\listing#1#2#3{{\sc #1}:\ {\it #2}, \ #3.}

\end{document}